\documentclass[12pt]{amsart}

\usepackage{am ssymb, stmaryrd, comment}
\usepackage{hyperref}
\usepackage[margin=1.0 in]{geometry}
\usepackage{color,soul}
\usepackage{multicol}
\usepackage[table]{xcolor}
\usepackage{enumerate}
\usepackage{mdframed}
\usepackage[T1]{fontenc}
\usepackage{tikz}
\usepackage{MnSymbol}
\usetikzlibrary{positioning}
\usepackage{multicol}
\usepackage{caption}
\usepackage{subcaption}
\usepackage{float}
\usepackage{cleveref}
\usepackage{tikz-cd}
\hypersetup{
    colorlinks=true,
    linkcolor=black,
    filecolor=black,
    citecolor= black,
    urlcolor=black}

\makeatletter
\def\bbordermatrix#1{\begingroup \m@th
  \@tempdima 4.75\p@
  \setbox\z@\vbox{%
    \def\cr{\crcr\noalign{\kern2\p@\global\let\cr\endline}}%
    \ialign{$##$\hfil\kern2\p@\kern\@tempdima&\thinspace\hfil$##$\hfil
      &&\quad\hfil$##$\hfil\crcr
      \omit\strut\hfil\crcr\noalign{\kern-\baselineskip}%
      #1\crcr\omit\strut\cr}}%
  \setbox\tw@\vbox{\unvcopy\z@\global\setbox\@ne\lastbox}%
  \setbox\tw@\hbox{\unhbox\@ne\unskip\global\setbox\@ne\lastbox}%
  \setbox\tw@\hbox{$\kern\wd\@ne\kern-\@tempdima\left[\kern-\wd\@ne
    \global\setbox\@ne\vbox{\box\@ne\kern2\p@}%
    \vcenter{\kern-\ht\@ne\unvbox\z@\kern-\baselineskip}\,\right]$}%
  \null\;\vbox{\kern\ht\@ne\box\tw@}\endgroup}
\makeatother

\def\VR{\kern-\arraycolsep\strut\vrule &\kern-\arraycolsep}
\def\vr{\kern-\arraycolsep & \kern-\arraycolsep}

\theoremstyle{plain}
\newtheorem{theorem}[subsection]{Theorem}
\newtheorem{prop}[subsection]{Proposition}

\newtheorem{lemma}[subsection]{Lemma}

\newtheorem{cor}[subsection]{Corollary}

\theoremstyle{definition}

\newtheorem{question}[subsubsection]{Question}

\newtheorem{example}[subsection]{Example}

\newtheorem{remark}[subsection]{Remark}

\newtheorem{construct}[subsection]{Construction}

\normalfont
\usepackage[T1]{fontenc}{}

\newcommand{\matlis}[1]{{#1}^{\vee}}
\newcommand{\Spec}{\text  {Spec}}
\newcommand{\Hom}{\text  {Hom}}
\newcommand{\End}{\text  {End}}
\newcommand{\frm}{\mathfrak{m}}

\title{Some algebras with trivial rings of differential operators}
\author{Alapan Mukhopadhyay}
\address{Institute of Mathematics, CAG,
EPFL SB MATH
MA A2 383 (Bâtiment MA), 
Station 8,
CH-1015 Lausanne,
Switzerland}
\email{alapan.mathematics@gmail.com}
\author{Karen E. Smith }
\address{Mathematics department, East Hall, 530 Church Street, University of Michigan, Ann Arbor- 48105, United States}
\email{kesmith@umich.edu}
\thanks{
This work is partially funded  by NSF DMS  \#2101075.}

\begin{document}

\begin{abstract} Let $k$ be an arbitrary field.  We construct examples of regular local $k$-algebras $R$ (of positive dimension) for which the ring of differential operators $D_k(R)$ is trivial in the sense that it contains {\it no} operators of positive  order. The examples are excellent in characteristic zero but not in positive characteristic. These rings can  be viewed as being non-singular but they are not simple as $D$-modules, laying to  rest speculation that $D$-simplicity might characterize a nice class of singularities in general. In prime characteristic, the construction also provides examples of {\it regular} local rings $R$ (with fraction field  a function field)  whose Frobenius push-forward $F_*^eR$ is {\it indecomposable} as  an $R$-module for all $e\in \mathbb N$. Along the way, we investigate hypotheses on a local ring $(R, m)$ under which  $D$-simplicity  for $R$ is equivalent to  $D$-simplicity  for its $m$-adic completion, and give examples of rings for which the differential operators do not behave well under completion. We also generalize a characterization of $D$-simplicity due to Jeffries in the $\mathbb N$-graded case: for a Noetherian local $k$-algebra $(R, m, k)$, $D$-simplicity of $R$ is equivalent to surjectivity of the natural map $D_k(R)\to D_k(R, k)$. 
\end{abstract}

\maketitle


\section{Introduction}
Let $R$ be a commutative algebra over a field $k$. The ring $D_k(R)$ of $k$-linear differential operators on $R$ is a natural non-commutative  extension of $R$  over which  $R$ may be viewed  as a module; see  \cite[IV \S16]{EGA}. For some time, the {\it simplicity}  of $R$ as a module over 
$D_k(R)$ (or "$D$-simplicity"  of $R$) has been thought to be related to  "mild" singularities of $\Spec R$. For example, in the setting of finite type $\mathbb C$-algebras, Levasseur and Stafford pointed out
an apparent connection between rational singularities and $D$-simplicity, and proved that many classical  rings of invariants are simple as $D$-modules \cite{LS89}. In prime characteristic,  for the local ring of a Frobenius split  point $x$ on a variety $X$,   $D$-simplicity is equivalent to   strong $F$-regularity\footnote{ Strong $F$-regularity can be viewed as a prime characteristic analog of  Kawamata log terminal singularities; see \cite{HaraWatanabe, SchwedeSmith}.}  at $x$ \cite{Sm95}.

This  note presents a cautionary tale: one cannot expect  the singularities of $\Spec (R)$ to be closely related to $D$-simplicity of $R$ much beyond the setting of  algebras essentially of finite type over $k$, or---in prime characteristic---much  beyond the $F$-finite\footnote{See \S \ref{primeChar} for the definition of $F$-finite.} case. Indeed, Theorem \ref{p: condition on derivative restricts differential operators} shows 
\begin{theorem}\label{Thm2} For any  field $k$, 
there exist regular local $k$-algebras $R$ of positive dimension which admit no differential operators of positive order---that is, the natural inclusion $R \hookrightarrow D_k(R)$ is an isomorphism.
\end{theorem}

 Such rings as in  Theorem  \ref{Thm2}, of course, are far from  $D$-simple:  any ideal of $R$ is a $D$-submodule. Furthermore, for these rings, the ring of differential operators is {obviously} finitely generated  over $R$. 
 While $D_k(R)$ is rarely finitely generated over $R$ (not even for the
$\mathbb C$-algebra  $\mathbb C[x, y, z]/(x^3+y^3+z^3)$ \cite{BGG}),  examples had suggested that for normal $\mathbb C$-algebras, finite generation of  $D_{\mathbb C}(R)$ over  $R$ might imply  $D$-simplicity of  $R$ \cite{LS89}.  Theorem  \ref{Thm2} destroys this naive hope as well, at least without additional  finiteness  assumptions. Examples of local $k$-algebras witnessing the above theorem are explicit and in particular do not admit any non-zero $k$-linear derivations; see \Cref{p: construction of g}, \Cref{c: DVRs with no derivations}.

For a prime characteristic ring $R$, the modules $F^e_*R$ obtained by viewing $R$ as a module over itself via the Frobenius map\footnote{See Section \ref{primeChar} for more careful introduction to the modules $F_*^eR$.} have played a central role in  prime characteristic singularity theory  since Kunz  proved that a Noetherian ring $R$ is regular if and only if $F^e_*R$ is  {\it flat} \cite{Kunz76}. In particular, when $R$ is the local ring of a point $x$  on a variety\footnote{i.e. integral finite type over an algebraically closed field.}  $X$, we see that $x$ is a smooth point if and only if $F^e_*R$ is a free $R$-module of rank $p^{e\dim  X}$;  more generally, $X$ is Frobenius split at $x$\footnote{Frobenius splitting  is closely related to log canonicity for points on a complex varieties;  see \cite{HaraWatanabe}.}  if and only if the $R$-module $F_*^eR$ admits a cyclic free summand. Hochster and Yao 
have proved, under some restrictions on the singularities and/or dimension,  that the $R$-modules $F_*^eR$  {\it {always}} decompose into a direct sum of (at least) two non-zero $R$-modules for $e\gg 0$, and asked\footnote{Schoutens also appears to have considered this issue quite generally \cite[footnote 7]{Schoutens}.}   whether this holds for {\it arbitrary} Noetherian $F$-finite rings of positive dimension \cite{Hochster, Yao, HochYau}.
Theorem \ref{Indecomp}  produces examples to  caution the reader that without the $F$-finite hypothesis,\footnote{A ring $R$ of characteristic $p$ is  $F$-finite if the $R$-module $F_*R$ is finitely generated.  An  example in \cite{DS18} suggested that  $F$-finiteness might be necessary for the indecomposability of $F^e_*R$, but no counterexample was known.} the answer to this question is negative, even for regular rings:

\begin{theorem}\label{Thm3}
There exist regular local rings of any prime characteristic and of positive dimension such that $F^e_*R$ is indecomposable for all $e\in \mathbb N$.
\end{theorem}

It is worth pointing out  that the examples witnessing Theorems \ref{Thm2}  and \ref{Thm3} are  local $k$-algebras whose ring of differential operators do not behave well under completion: $\widehat R \otimes_R D_k(R)\not\cong D_{k}(\widehat R)$ (Corollary \ref{pathological}). On the other hand, 
for local $k$-algebras $(R, m)$  that are either essentially of finite type over $k$ or   $F$-finite,   $\widehat R \otimes_R D_k(R)\cong D_{k}(\widehat R)$ (Corollary \ref{p: differential operator of a completion} and Proposition \ref{p: for an F-finite ring differential operators of completion is the scalar extension of differential operators}). In these cases, assuming furthermore that $R/m\cong k$,  
we show that $D$-simplicity is equivalent to $D$-simplicity after completion (Corollary \ref{c: largest D ideal and completion}). This follows from  Theorem \ref{MainThm}, which   
generalizes a characterization of $D$-simplicity of Jeffries  in the graded case  to arbitrary Noetherian local rings with a coefficient field.

The rings we construct for Theorems \ref{Thm2} and \ref{Thm3} are fairly nice discrete valuation rings  (with a generalization to dimension two) refining a construction from \cite{DS18}.  Their fraction fields are finitely generated over $k$, so they are {\it generically}  essentially of finite type over $k$.   There is a natural way to generalize the construction to higher dimension which produces rings $R$  witnessing the behavior in Theorems  \ref{Thm2} and \ref{Thm3} {\it when  $R$ is  Noetherian. }  The Noetherianity is known, however, only in dimension two \cite{Valabrega}; see Remark \ref{higherdim} and Question \ref{Q}.

  In prime characteristic, our work adds to the growing body of evidence that   {\it $F$-finite} schemes are the natural setting in which to study $F$-singularities: $F$-finite schemes appear to behave very much like finite type $k$-schemes, while expected properties tend to fail for non-$F$-finite schemes.
 For example, the aforementioned results connecting  the structure of $F_*^eR$ to classes of singularities all fundamentally use the fact that the local 
 rings of varieties are  $F$-finite. 
  In a slightly different direction, it is known that for $F$-finite rings $R$,   the module $F_*^eR$, for $e\gg 0$,  strongly generates the bounded derived category of Noetherian $R$-modules in the sense of Bondal-Van den Bergh;  see \cite{BLIMP}.
  Strikingly, $F$-finite schemes have a {\it canonically defined} dualizing complex, as well \cite{Bhattetaldualizingcomplex}.

\section{Differential operators and completion}

 We recall the definitions from \cite[IV \S 16]{EGA}. Fix a commutative ground ring  $A$ and a commutative  
 $A$-algebra $S$.
  Let $N$ and $M$ be  $S$-modules.  The $A$-linear 
differential operators from $N$ to $M$, denoted $D_{S/A}(N, M), $   form an $A$-submodule of $\Hom_A(N, M)$, with a natural
 left (respectively, right) $S$-module structure induced from  $S$-module structure on $M$ (respectively, $N$). 
The $S$-bimodule $D_{S/A}(N,M)$ of $A$-linear differential operators from $N$ to $M$ is defined inductively as follows:
First, the differential operators of order zero are simply the $S$-linear maps 
$$D^0_{S/A}(N,M)= \text{Hom}_S(N,M).$$
 For $n \geq 1$, we say that $\delta\in \Hom_A(N, M)$ is in 
 $D^n_{S/A}(N,M)$, the $A$-module of differential operators of order at most $n$,  if and only if  for all $s \in S$, the commutator $[\delta, s]$ is in $D^{n-1}_{S/A}(N,M)$. 
 Then $$D_{S/A}(N,M)= \underset{n \in \mathbb{N}}{\bigcup}D^n_{S/A}(N,M) \subseteq \text{Hom}_A(N,M)$$
is the full $S$-bimodule of differential operators from $N$ to $M$.

The $A$-module  $D_{S/A}(M,M)$ has a natural (non-commutative) ring structure under composition. In particular, the ring $D_{S/A}(S,S)$, usually denoted by $D_A(S)$, is called \textit{the ring of $A$-linear differential operators of $S$}.  Note that $D_{S/A}(N,M)$ has a natural structure of a right $D_{S/A}(N)$-module as well as a left $D_{S/A}(M)$-module, given by pre and post composition with elements of  $D_{S/A}(N,M)$.   We sometimes drop the $S$  from the notation in $D_{S/A}(N,M)$ when the ring is clear from the context.

\subsection{Functoriality}\label{functorial} Differential operators are  functorial in both $N$ and $M$: fixing $M$, any $S$-module map $N\to N'$ induces an $S$-bimodule (or left $D_{S/A}(M)$-module) map $D_{S/A}(N', M)\to D_{S/A}(N, M)$; likewise, fixing $N$, any $S$-module map $M\to M'$ induces an $S$-bimodule (or right  $D_{S/A}(N)$-module)  map $D_{S/A}(N, M)\to D_{S/A}(N, M')$.

\subsection{Localization}\label{localize}  Differential operators from $M$ to $N$   extend naturally without fuss  to localization at any multiplicative system $W$, similar to the "quotient rule" in calculus. Indeed,
 a zero-order operator $M\to N$ (that is, an $S$-linear map) extends to an $S$-linear map  $W^{-1}M\to W^{-1}N$ by tensoring with $W^{-1}S$ over $S$. 
So, assuming inductively that operators in $ D_{S/A}^n(M, N)$ extend uniquely   to operators in  $ D_{S/A}^n(W^{-1}M, W^{-1}N)$, we can extend $\phi\in  D_{S/A}^{n+1}(M, N)$ to the localization $D_{S/A}^{n+1}(W^{-1}M, W^{-1}N)$ using the formula\footnote{The formula (\ref{localization}) in fact gives a differential operator in $D^{n+1}_{W^{-1}S/A}(W^{-1}M, W^{-1}N)$.}
\begin{equation}\label{localization}
\phi(\frac{m}{s}) = \frac{\phi(m)-[\phi, s](\frac{m}{s})}{s}  \in W^{-1}N
\end{equation}
 for each $\frac{m}{s}\in W^{-1}M$.

\subsection{Universal Differential operators}\label{Princ}
Alternatively,  we can define $D_{S/A}(N, M)$ as follows.  Let  $\Delta_{S/A}$ be the kernel of the multiplication map $S \otimes_A S \rightarrow S$. For each $n \in \mathbb{N}$,  the module of $n$-th {\bf principal parts} for $N$ is defined as
$$P^n_{S/A}(N)= \frac{S \otimes_A S}{\Delta_{S/A}^ {n+1}}\otimes_S N;$$
note that $P^n_{S/A}(N)$ 
has a natural  $S$-bimodule structure from the left and from the right factors.
Thinking of  $P^n_{S/A}(N)$ as a {\it left} $S$-module, consider  the natural $A$-module map
 \begin{equation}\label{universal}
 N \xrightarrow{\iota_n}
  P^n_{S/A}(N)\,\,\,\,\,\,\,\,\, {\text{ sending }} \,\,\,\,\,\,\,\,\,\, x\mapsto (1\otimes_A 1 \,\, {\text{mod}} \,\, \Delta_{S/A}^{n+1})\,   \otimes_S x. 
 \end{equation}
 The map $\iota_n$ is the {\bf universal $A$-linear   differential operator of order $n$} on $N$: 
  Composition with any left $S$-module map $P^n_{S/A}(N) \to M$  produces an 
  element of $\text{Hom}_A(N,M)$ which is a differential operator of order at most  $n$, and conversely, 
  every differential operator in $D^n_{S/A}(N, M)$ factors through (\ref{universal}). Put differently, considering $P^n_{S/A}(N)$ as a {\it left}  $S$-module,
   there is a natural isomorphism  of left\footnote{The map (\ref{e: differential operators are dual of module of principal parts}) is also an isomorphism of {\it right} $S$ modules, where the right structure on  $\text{Hom}_S(P^n_{S/A}(N),M)$ comes from the right $S$-module structure on $P^n_{S/A}(N)$,  but we won't make use of this structure.}
   $S$-modules
\begin{equation}\label{e: differential operators are dual of module of principal parts}
\text{Hom}_S(P^n_{S/A}(N),M) \rightarrow D^n_{S/A}(N,M) \, \, \, \,\,\, \textup{sending}\,\,\,\,\,\, \phi \mapsto \phi \circ \iota_n.
\end{equation}
\\

\subsection{Continuity of Differential Operators} For any ideal $J$ in  $S$, differential operators $N\to M$ are continuous in the $J$-adic topology.\footnote{The reader may consult  \cite[Chapter 10]{AM} for basic definitions and facts about the $J$-adic topology and completion.} This has the following consequence:

\begin{prop}\label{p: differential operator of completion without any noetherian assumption}
    Let $S$ be an $A$-algebra and $J$ an ideal of $S$. Let $M$ and $N$ be $S$-modules, with $M$ $J$-adically complete. Then each $A$-linear differential operator $\phi:N\to M$ extends  {\it uniquely} to an $A$-linear differential operator $\phi:\widehat N^J\to M$ of the same order.
    Put differently, for each $n\in \mathbb N$,  the natural map 
    \[
    D^n_{\widehat{S}^J/A}({\widehat N}^J, M) \rightarrow D^n_{S/A}(N, M)
    \] induced by the $J$-adic completion map $N \rightarrow \widehat {N}^J$ is a bijection  respecting the $S$-bimodule structure.
\end{prop}

To prove Proposition \ref{p: differential operator of completion without any noetherian assumption}, we need the following well-known lemma:

\begin{lemma}\label{l: continuity of differential operators}\cite[Prop 1.4.6]{Yek92}
 Fix an $A$-algebra $S$ and an  ideal $J\subseteq S$. Let $M$ and $N$ be $S$-modules, and let  $\phi: N\to M$ be a differential operator of order at most $r$. 
 For all $n\in \mathbb N$,  $\phi(J^{n}N) \subseteq J^{n-r}M$  (where we interpret any non-positive power of $J$ as the unit ideal $S$).\end{lemma}


\begin{proof}
We induce on $n+r$.   Without loss of generality, we can assume $n \geq r$. When $r=0$, the operator $\phi$ is $S$-linear; thus $\phi(J^{n}N) \subseteq J^{n}M$.  Assume that $r \geq 1$.   Observe that the abelian group $J^nN$ is  generated by elements of the form
$ss'x$, where $s'$ and $s$ are in $J^{n-1}$ and $J$ respectively, and $x\in N$. Note that
 $$\phi(ss'x)= s\phi(s'x)+ [\phi, s](s'x) $$
and $[\phi, s] \in D^ {r-1}_A(N,M)$. Both $s \phi(s'x)$ and $[\phi, s](s'x)$ are in $J^{n-r}M$  by our induction hypothesis, hence also  $\phi(ss'x)\in J^ {n-r}M$.
\end{proof}

\medskip

Before continuing with  the proof of  Proposition \ref{p: differential operator of completion without any noetherian assumption},  we point out the following useful application of Lemma \ref{l: continuity of differential operators}:

\begin{example}\label{ContMapsAreDiffOps} Fix a field $k$.
Let $(R,m)$ be a  local $k$-algebra such that the natural composition $k\hookrightarrow R \to R/m$ is an isomorphism.
Then  as subsets of $\Hom_k(R, R/m)$,  \[D_{R/k}(R, R/m) = \Hom_k^{cts}(R, R/m),\] the $k$-vector space of {\it  all} $m$-adically continuous  $k$-linear maps $R\to R/m$.
\end{example}

\begin{proof}[Proof of Example \ref{ContMapsAreDiffOps}]  Lemma \ref{l: continuity of differential operators} ensures that 
 $$D_{R/k}(R, R/m) \subseteq \Hom_k^{cts}(R, R/m).$$ On  the other hand, if $\phi:R\to R/m$ is continuous in the $m$-adic topology, then there exists $n\in \mathbb N$ such that  $\phi(m^{n}) = 0$. We prove by induction on $n$  that this means $\phi$ is a differential operator of  order at most $n-1$. 
 
 First: by hypothesis, the natural surjection $R\to R/m$ is a splitting of  $k\hookrightarrow R$, so that as a $k$-vector space, 
 $R\cong k \oplus m$. So we can write any $r\in R$ as $r=\lambda +x$ where $\lambda\in k$ and $x\in m$.

 Now to check the base case of  our induction, assume that $\phi(m)=0$. In this case, taking arbitrary  $r\in R$, 
  and writing 
 \[
 r = \lambda  + x,\,\,\,\,\,\,\,\,\,\,\,\,{\text{and}} \,\,\,\,\,\,\,\,\,\,\,\, s = \lambda'  + x'
 \]
 where $\lambda, \lambda'\in k$ and $x, x'\in m$, we have
 \[
 s\phi(r) = (\lambda' +x') (\phi(\lambda) + \phi(x)) = \lambda' \phi(\lambda) = \phi(\lambda'\lambda),\]
 which is $\phi( \lambda'\lambda  +x' \lambda + \lambda' x  +x'x)$, or $\phi(sr)$.
  That is, $\phi$ is $R$-linear and hence a differential operator of order zero.
 
 Now assume that $\phi(m^{n+1}) = 0$. We must show $\phi\in D_k^n(R)$. For this, we  fix  $r\in R$, and consider the operator $[\phi, r]$. Writing $r=\lambda +x$ where $\lambda\in k$ and $x\in m$, we see 
 \[
[\phi, r] = [\phi, \lambda + x ] =   [\phi, \lambda] + [\phi,  x ] = 0 +   \phi \circ x - \overline{x}\circ \phi  =  \phi \circ x,
 \]
 where $\overline x$ denotes the class of $x$ in $R/m$.
Since $x\in m$ and $\phi(m^{n+1})=0$,  we see that $\phi \circ x$ kills $m^{n}$. By  induction, we conclude that $[\phi, r]\in D^{n-1}(R).$ Since $r\in R$ was arbitrary, we conclude that $\phi\in D^n(R)$.
 \end{proof}

\begin{proof}[Proof of Proposition \ref{p: differential operator of completion without any noetherian assumption}] 
It suffices to prove\footnote{To reduce notational clutter, we drop  $J$ in the notation  $\widehat S^J$ and $\widehat N^J$  in this proof.} the first sentence: any $A$ differential operator of $S$-modules  $N\to M$ extends   uniquely to an $A$-linear differential operator of  $\widehat S$-modules $\widehat N\to M$.  

First note that  by definition,  because $M$ is $J$-adically complete,
the natural map $M\to \varprojlim M/J^tM$ is an isomorphism. In particular, $\bigcap_{t\in \mathbb N} J^tM  = 0. $ 
By  Lemma \ref{l: continuity of differential operators}, 
 the image of $\underset{t \in \mathbb{N}}{\bigcap}J^tN$ under any differential operator $N\to M$  is contained in 
$\underset{t \in \mathbb{N}}{\bigcap}J^{t}M$, so  the image of $\underset{t \in \mathbb{N}}{\bigcap}J^tN$ under any differential operator is zero. 
 Hence, for every $r\in \mathbb N$,  the quotient map $$N \longrightarrow N/\underset{t \in \mathbb{N}}{\bigcap}J^tN$$ induces a natural bijection  $ D^r_A( N/\underset{t \in \mathbb{N}}{\bigcap}J^tN,M) \to D^r_A(N,M)$. Thus replacing $N$ by $ N/ \underset{t \in \mathbb{N}}{\bigcap}J^tN$, we may assume without loss of generality that the completion map $N \rightarrow \widehat {N}$ is injective.\\

Now fix a differential operator $\phi: N\to M$ of order at most $n$.  We wish to extend  $\phi$  to a differential operator $\widehat\phi:  \widehat N\to M.$  For each $x \in \widehat{N}$,   choose a Cauchy sequence $(x_m)_m$ in $N$  converging to $x$ in $ \widehat{N}$. The sequence  $\left(\phi(x_m)\right)_m$ in $M$  is Cauchy by \Cref{l: continuity of differential operators}, so has a uniquely defined limit  $\underset{m \to \infty}{\lim} \phi(x_m)$ in the $J$-adically complete module $M$. Now define $\widehat \phi$ by
\[
\widehat \phi: 
\widehat{N} \to M \,\,\,\,\,\,\,\,\,\,\, \,\,\,\,\,\,\,\,\,\,\,  x\mapsto   \underset{m \to \infty}{\lim} \phi(x_m).
\]
The map $\widehat \phi$  is clearly $A$-linear and restricts to $\phi$ on  $N$. We  must show that $\widehat \phi$  is a differential operator (of $\widehat S$-modules)  of order at most $n$. That is, we must  show that $\widehat \phi$ factors through the universal differential operator $\widehat N\to P^n_{\widehat S/A}(\widehat N)$ (see Paragraph \ref{Princ}). For this,
 it suffices to check that the left  $\widehat S$-module map 
\begin{equation}\label{map3}
\begin{aligned}
\widehat S \otimes_A \widehat S \otimes_{\widehat S} \widehat N &\xrightarrow{\Phi} M\\
y \otimes_A 1 \otimes_{\widehat  S}  x & \mapsto     y \widehat{\phi}(x)
\end{aligned}
\end{equation}
 satisfies $\Phi(\Delta^{n+1}_{\widehat S/A}  \otimes_{\widehat S} N)=0.$ 

To this end, for each $s \in \widehat S$, set $d(s)= s \otimes 1- 1 \otimes s \in \widehat S \otimes_A \widehat S$, and  
recall that the collection of all such $d(s)$ generate 
$\Delta_{\widehat S /A}$ as an ideal of  $\widehat S\otimes_A \widehat S$. Since $\Phi$ is left $\widehat S$-linear, then, it suffices to show that each "monomial" of the form 
\[
    d(s_1)d(s_2)\cdots d(s_{n+1}) \otimes_{\widehat  S} x \in \widehat S \otimes_A \widehat S \otimes_{\widehat S^J} \widehat N^J
 \]
 is taken to zero by $\Phi$.  
 We accomplish  this by  showing that for every $t\in \mathbb N$, 
 \begin{equation}\label{need}
  \Phi(d(s_1)d(s_2)\cdots d(s_{n+1}) \otimes_{\widehat  S} x) \in J^{t-n} M, 
 \end{equation}
 from whence it follows that $\Phi(d(s_1)d(s_2)\cdots d(s_{n+1}) \otimes_{\widehat  S} x) \in \bigcap_{t\in \mathbb N} J^{t} M = 0$.
 
  To prove (\ref{need}),  we write each  $s_j \in \widehat S$ (respectively, $x\in \widehat N$) 
  as the limit of a Cauchy sequence $(s_{jt}\, | \,  t \in \mathbb N)$ in $S$ (respectively,  $ (x_{t}\, | \,  t \in \mathbb N) $ in $N$) such that 
    \[
  s_j - s_{jt} \in J^t\widehat S\,\,\,\,\,\,\,\,\, {\text{ and }} \,\,\,\,\,\,\,\,\,  x-x_t \in J^t\widehat N
  \]
for all  $t\in\mathbb N$.  Since $d$ is additive, one easily checks that for all $t\in \mathbb N$, 
\[
\begin{aligned}
    d(s_1)d(s_2)\cdots d(s_{n+1})  &\in  d(s_{1t})d(s_{2t})\cdots d(s_{n+1\ t }) +  J^t\widehat S \otimes_A \widehat S + \widehat S \otimes_A J^{t}\widehat S, \,\,\,\,\, {\text{and}} \\
     d(s_1)d(s_2)\cdots d(s_{n+1})  \otimes x  &\in  d(s_{1t})d(s_{2t})\cdots d(s_{n+1\ t }) \otimes x_t  +
      J^t\widehat S \otimes_A \widehat S \otimes_{\widehat S} \widehat N + \widehat S \otimes_A J^{t}\widehat S    \otimes_{\widehat S}  \widehat N.
      \end{aligned}
\]
Now applying the  left $\widehat S$-linear map $\Phi$ to the monomial $d(s_1)d(s_2)\cdots d(s_{n+1})  \otimes x$, we have
\[
\begin{aligned}
\Phi( d(s_1)d(s_2)&\cdots d(s_{n+1})  \otimes x )  \in \\ 
 &\underbrace{\Phi(d(s_{1t})d(s_{2t})\cdots d(s_{n+1\ t }) \otimes x_t )}_{{\text{in } \Delta^{n+1}_{S/A} \otimes_S N}} + \
     \Phi(J^t\widehat S \otimes_A \widehat S \otimes_{\widehat S} \widehat N)  +  \Phi(\widehat S \otimes_A J^{t}\widehat S    \otimes_{\widehat S}  \widehat N)\\
     &\,\,\,\,\,\,  =  0 + J^t \widehat \phi (\widehat N)  +  \widehat \phi(J^{t}  \widehat N) \subseteq J^{t-n} M.
     \end{aligned}
\]
Here, the first term is  zero  because $\phi:N\to M$ is a differential operator (of $S$-modules)  of order at most $n$.  And the last inclusion follows because $\widehat \phi (J^tN)\subseteq J^{t-n}M$ for all $t> 0$. Indeed, any $y\in J^t \widehat N$ is a limit of a Cauchy sequence $(y_m)_m$ where $y_m \in J^tN$ and the $J$-adic continuity of $\phi$ implies $\phi(y_m) \in J^{t-n}M$, see \Cref{l: continuity of differential operators}. Since $J^{t-n}M$ is complete in the subspace topology $\phi(y) \in J^{t-n}M$; see \cite[cor 10.3]{AM}.
\end{proof}

\begin{cor}\label{p: differential operator of a completion} 
Let $S$ be an  algebra essentially finite type over a Noetherian ring $A$, and let $J$ be an ideal of $S$.  Then there is a natural left $\widehat S^J$-module isomorphism
\[
  \widehat S^J \otimes_S D^n_A(S) \rightarrow D^n_A(\widehat S^J),  
\]
where the tensor product is formed using the left $S$-module structure on $D^n_A(S)$. Moreover, $\widehat S^J\otimes_S D_A(S)\cong D_A(\widehat S^J)$ as rings.
\end{cor}

\begin{proof}
 The inclusion of targets $S \subseteq \widehat S^J$ yields a natural map $D^n_A(S,S) \rightarrow D^n_A(S, \widehat S^J)$. Following with
the natural isomorphism $D^n_A(S, \widehat S^J) \rightarrow D^n_A(\widehat S^J, \widehat S^J)$ of  Proposition
 \ref{p: differential operator of completion without any noetherian assumption} 
produces a  map
 $D^n_A(S)$ to $D^n_A(\widehat S^J)$ of left $S$-modules. So we have a natural map of $\widehat S^J$-modules
 \begin{equation}\label{e: the isomorphism of differential operators of the completion}
  \widehat S^J \otimes_S D^n_A(S) \rightarrow D^n_A(\widehat S^J), 
\end{equation}
or equivalently,
\begin{equation}\label{e: using the finite generation of module of pricipal parts}
\widehat S^J \otimes_S \text{Hom}_S(P^n_A(S),S) \rightarrow \text{Hom}_{\widehat S^J}(\widehat S^J \otimes_S P^n_A(S),\widehat S^J).   
\end{equation}
Because $S$ is essentially of finite type over the Noetherian ring  $A$, the $S$-module $P^n_A(S)$ is finitely presented over $S$, and so the map (\ref{e: using the finite generation of module of pricipal parts}), and hence (\ref{e: the isomorphism of differential operators of the completion})  is an isomorphism. Finally, it is straightforward to verify that the ring structure in the union is preserved as well. 
  \end{proof}

\subsection{$F$-finite rings}\label{primeChar}
Corollary \ref{p: differential operator of a completion}  holds for Noetherian rings of prime characteristic 
without the finite type assumption, provided that the ring is {\it  $F$-finite.} Recall that 
every commutative ring of prime characteristic $p>0$ is endowed with  natural  ring homomorphisms 
$$
R \overset{F^e}\to R \,\,\,\,\,\,\,\,\, \,\,\,\,\,\,\,\,\, \,\,\,\,\,\,\,\,\, r\mapsto r^{p^e}
$$
called  the Frobenius maps, which  induce   $R$-module structures on the target copy of $R$ via restriction of scalars. We denote these $R$-modules by $F_*^eR$, and declare $R$ to be {\bf $F$-finite} if $F_*^eR$ is finitely generated for some (equivalently, all) $e\in \mathbb N$.

There is a well-known\footnote{This is clearly stated in \cite[Thm 1.4.9]{Yek92}, with partial attributions to \cite[Lem 3.3]{Chase} and \cite[Proof of Thm 1]{Wod}. See 
 \cite[\S 2.5]{SVdB} for a straightforward proof; note that although $A$ is assumed a field in \cite{SVdB}, it is never used.}  description
 of rings of differential operators for $F$-finite rings:
 \begin{prop}\label{p: alternate description of the ring of differential operators of a prime characteristic ring} Let $A$ be a perfect ring of prime characteristic $p$ and let $S$ be an $A$-algebra.
\begin{enumerate}
    \item The ring of differential operators $D_A(S)$ is contained in $\underset{n \in \mathbb{N}}{\bigcup}\text{End}_S(F^n_*S)$.
    \item If $S$ is  $F$-finite, then  $D_A(S)= \underset{n \in \mathbb{N}}{\bigcup}\text{End}_S(F^n_*S)$.
\end{enumerate}
\end{prop}

Using the description of Proposition \ref{p: alternate description of the ring of differential operators of a prime characteristic ring}, we see that differential operators commute with completion for $F$-finite rings: 

\begin{prop}\label{p: for an F-finite ring differential operators of completion is the scalar extension of differential operators}
 Let $A$ be a  perfect\footnote{A prime characteristic ring is {\bf perfect} if $A^p=A$. Bhargav Bhatt pointed out that if a connected Noetherian ring contains a  prime characteristic perfect ring which is not necessarily Noetherian, the perfect ring must be a domain and its fraction field is contained in the ambient Noethrian ring. So a prime characteristic perfect subring of a Noetherian ring is contained in a finite product of perfect fields, where the fields are contained in the Noetherian ring.} ring, and let $S$ be a Noetherian
 $F$-finite $A$-algebra. Then completing at any ideal $J$ of $S$, 
the natural map $$\widehat S^J \otimes_S D_A(S) \rightarrow D_A(\widehat S^J)$$ is a ring isomorphism. 
 \end{prop}
 
 \begin{proof}  Under the hypotheses of Proposition  \ref{p: for an F-finite ring differential operators of completion is the scalar extension of differential operators}, the maps $$\widehat S^J\otimes_S \text{End}_S(F^n_*S)\to
  \text{End}_{\widehat S^J}(\widehat S^J \otimes_SF^n_*S)   \to   \text{End}_{\widehat S^J}(F^n_*\widehat S^J) $$
  are isomorphisms for all $n$: the first arrow is an isomorphism because $F^n_*S$ is finitely presented for every $n\in \mathbb N$  \cite[Tag 00F4]{Stacksproject}), and the second by, {\it e.g.}  
 \cite[Cor 1.15]{SchwedeSmith}. So since $\widehat S^J$ is $F$-finite  \cite[Prop 1.21(e)]{SchwedeSmith}, taking the limit over all $n\in \mathbb N$, it follows that
$\widehat S^J\otimes_S D_k(S)\cong  D_k(\widehat S^J)$.
 \end{proof}

\begin{remark} The behavior of 
differential operators under completion is described in \cite[Thm p273]{Musson}, \cite[Cor p13]{Ishibashi} and  \cite[page 45]{Lyubeznik} in various more restrictive settings.  \end{remark}

\section{A characterization of $D$-simplicity}
In this section, we prove the following  characterization of $D$-simplicity  for  a local Noetherian  $k$-algebra  "with a coefficient field".  This includes  the local rings of $k$-rational   points on finite type $k$-schemes, as well as complete local Noetherian rings (by the Cohen Structure Theorem \cite{Cohen}). 

\begin{theorem}\label{MainThm} Fix a field $k$, and 
let $(R, m)$ be a Noetherian local $k$-algebra such that the natural composition $k\hookrightarrow R \to R/m$ is an isomorphism. Then  $R$ is simple as a $D_k(R)$-module if and only if 
 the natural map 
\begin{equation}\label{NatMap1}
D_k(R) \to D_k(R, R/m)
\end{equation}
induced by the residue map  $R\to R/m$ is surjective.
\end{theorem}

Theorem \ref{MainThm}  generalizes a result of Jeffries, who proved a graded version,   though his proof also works in the complete case; see  
  \cite[Proposition 3.2]{Jeff21}. Our proof holds  more generally and provides  additional insight into the cokernel of (\ref{NatMap1}), which we will describe  in  terms of  following obvious obstruction to $D$-simplicity. Similar obstructions to $D$-simplicity have been considered more generally in \cite[section 3]{BJB19}.
\begin{prop}\label{l: larget D-ideal}  Let $(R,m)$  be a local  $k$-algebra. 
Then
the ideal 
\begin{equation}\label{M}
\mathcal M_R= \{x \in R \, | \, \delta(x) \in m, \, \forall \, \delta \in D_k(R)\}
\end{equation}
 is the unique  
maximal  proper $D_k(R)$-submodule of $R$.
\end{prop}
\begin{proof}
    The set $\mathcal M_R$ is proper  because $1 \in R$ is mapped to $1$ by the identity operator, thus $1 \notin \mathcal M_R$.  The set $\mathcal M_R$ is a $D_k(R)$-submodule of $R$ because for arbitrary  $x \in \mathcal M_R$ and $\phi \in D_k(R)$, we have 
    $\delta(\phi(x))= (\delta \circ \phi) (x) \in m$ for any $ \delta \in D_k(R)$. In particular, $\mathcal M_R$ is an ideal of $R$. This ideal is the unique largest $D_k(R)$-submodule of $R$, because every  proper $D_k(R)$-submodule $N \subseteq R$ must satisfy $\delta(N) \subseteq N \subseteq m$ for all $\delta \in D_K(R)$: otherwise, $N$ would contain a unit and hence be equal to $R$. 
\end{proof}

In light of Proposition \ref{l: larget D-ideal}, Theorem \ref{MainThm} follows immediately from the faithful exactness of Matlis duality and the following  result:

\begin{theorem}\label{c: main theorem}
Fix a field $k$. Let $(R,m,k)$ be 
 a Noetherian local $k$-algebra  such that the natural composition $k \hookrightarrow R \rightarrow \frac{R}{m}$ is an isomorphism.
The Matlis dual 
of the natural (right) $R$-module map 
\begin{equation}\label{NatMap}
D_k(R) \to D_k(R, R/m)
\end{equation} 
is a map  \[
\widehat  R \to \matlis{D_k(R)}
\]
whose kernel is  $\mathcal M_R\widehat R$,  where $\mathcal M_R$ is the maximal proper $D$-submodule of $R$ as defined in (\ref{M}).
\end{theorem}

The proof of Theorem \ref{c: main theorem} will follow a digression on Matlis duality.

\begin{remark}
Our statement and proof  of Theorem \ref{c: main theorem} and  Proposition \ref{l: larget D-ideal} adapt to the graded category with appropriate minor modifications. See Remark \ref{graded}.
\end{remark}

\medskip

\subsection{Injective Hulls and Matlis duality}\label{InjHullMatDual} Matlis duality is a well-known tool in commutative algebra;  see, for example,  \cite[\S18]{Matsumura}. Here, we summarize the key points needed to prove Theorem \ref{c: main  theorem}.

 Let $(R,m)$ be a Noetherian local $k$-algebra such that the composition $k\hookrightarrow R \to R/m$ is finite. 
Then the $R$-module  
 \begin{equation}\label{injHull}
E := \text{Hom}^{cts}_k(R,k)\cong \underset{n}{\varinjlim} \,\text{Hom}_k(\frac{R}{m^n},k),
 \end{equation}
 of $m$-adically continuous maps in  $\text{Hom}_k(R,k)$ is an {\it injective hull of the residue field of $R$}; see {\it e.g.}  \cite[Rmk 3.2]{Smith18}.  
 In  particular, under the hypothesis of Theorem \ref{c: main  theorem}, Example \ref{ContMapsAreDiffOps}  implies that
  \begin{equation}\label{E=D}
 E  :=   \text{Hom}^{cts}_k(R,k) \cong D_k(R, R/m)
 \end{equation}
 as submodules 
 of  $\text{Hom}_k(R,k)$.
 The $R$-module $E $ evidently has a natural $\widehat R$-module structure;  as such, $E $ is also  an injective hull of the residue field of $\widehat R$.
 
The  \textit{Matlis dual} functor is the exact, faithful contravariant functor on the category of $R$-modules
\[
M \,\, \mapsto \,\,\, \matlis{M}:=\text{Hom}_R(M, E).
\]
Because $E$ is naturally an $\widehat R$-module, Matlis duality takes $R$-modules to $\widehat  R$-modules,  and the functors 
\begin{equation}\label{9}
\Hom_{R}(-, E) \,\,\,\,\,\,\,\,\,\, {\text{and}} \,\,\,\,\,\,\,\,\, \Hom_{\widehat  R}(\widehat R \otimes_R -, E)
\end{equation}
are equivalent. Moreover,  $E^\vee \cong \widehat R$, and for every Noetherian $R$-module $M$, 
the natural  map $M \to (M^\vee)^{\vee}$ agrees with completion $M \to\widehat R\otimes_{R}M$.  In our setting, even when $M$ is Noetherian, we can identify $\matlis{M}:=\text{Hom}_R(M, E)$ with 
$ \text{Hom}^{cts}_k(M, k),$ the  submodule of maps in $ \text{Hom}_k(M, k)$ continuous in the $m$-adic topology; see \cite[Prop 3.1]{Jeff21}.

\begin{remark}
In fact, the Matlis dual functor can be viewed as a functor from left $D_k(R)$-modules to right  $D_k(R)$-modules (and vice-versa). This point of view is developed by Switala in \cite[\S 4]{Swi}  and is used by Jeffries in his proof of Theorem \ref{MainThm} in the graded case. We will not make use of it in this paper.  
\end{remark}

The proof of  Theorem \ref{c: main theorem} requires the following technical lemma:
\begin{lemma}\label{l: kernel of completion}
Let $J$ be an ideal in a Noetherian ring $S$. Let $N$ be a finitely generated $S$-module and $M$ be a module over the completion $\widehat S^J$.  Let $\psi: N \rightarrow M$ be an $S$-linear map. Then the kernel of the map obtained from $\psi$ by extending scalars:
\[N \otimes_S \widehat S^J \rightarrow M\]
is isomorphic to $\textup{ker}(\psi)\otimes_S \widehat S^J$.
\end{lemma}
\begin{proof}
Let $N'$ be the $\hat S^J$-submodule of $M$ spanned by $\psi(N)$.
We claim that the map \[\theta: \psi(N) \otimes_S \widehat S^J \rightarrow N',\] obtained from the inclusion $\psi(N) \hookrightarrow N'$ by extension of scalars, is an isomorphism.

To establish the claim, note that the $S$-submodule $\psi(N)$ is dense in $N'$ in the $J$-adic topology. So the $S$-linear map 
$$\psi(N) \rightarrow \psi(N) \otimes_S \widehat S^J \, \, \, \,\,\, \textup{sending}\,\,\,\,\,\, x \mapsto x \otimes 1,$$
uniquely extends to an $\widehat S^J$-linear map $\theta': N' \rightarrow \psi(N) \otimes_S \widehat S^J$. Indeed, given $x \in N'$, choose a sequence of elements $(x_n)_n$ in $\psi(N)$ converging to $x$, then define $\theta'(x)$ to be the limit of $(x_n \otimes 1)_n$ in  $\psi(N) \otimes_S \widehat S^J$. The limit of $(x_n \otimes 1)_n$ exists and is independent of the choice of $(x_n)_n$ as $\psi(N) \otimes_S \widehat S^J$ is complete and separated in the $J$-adic topology, thanks to the finite generation hypothesis on $N$. The $S$-linear map $\theta' \circ \theta$ restricted to $\psi(N) \subseteq \psi(N)\otimes \widehat S^J$ is the identity map. Now $\psi(N)\otimes \widehat S^J$ is the $J$-adic completion of its submodule $\psi(N)$. The endomorphism $\theta' \circ \theta$ is the unique continuous extension of $\theta' \circ \theta_{|\psi(N)}$. So $\theta' \circ \theta$ is the identity map. Therefore $\theta$ is injective. As $\theta$ is clearly surjective, the claim follows.

Now, the map $N \otimes_S \widehat S^J \rightarrow M$ in the statement of the lemma factors as
\[N \otimes_S \widehat S^J \rightarrow \psi(N)\otimes_S \widehat S^J \cong N' \hookrightarrow M.\]
So the desired kernel is the kernel of the map $N \otimes_S \widehat S^J \rightarrow \psi(N)\otimes_S \widehat S^J$ induced by $\psi$. By the flatnss of $\widehat S^J$, 
this means that the kernel is isomorphic to $\textup{ker}(\psi)\otimes \widehat S^J$.
\end{proof}

\begin{proof}[Proof of Theorem \ref{c: main theorem}]
     Consider the map $D_k(R, R) \to D_k(R, R/m)=E$ induced by the natural surjection $\pi:R\to R/m$ (see \S \ref{functorial}). Applying the Matlis dual functor, to this map, we get an $\widehat R$-linear map
     \begin{equation}
     \phi: \widehat R \cong \textup{Hom}_R(E,E) \rightarrow \textup{Hom}_R(D_k(R),E) \ .
     \end{equation}
We identify $E$ with $\textup{Hom}_k^{cts}(R,k)$. For an element $r$ in $R$, $\phi(r)$ maps a differential operator $\delta$ to $\pi
\delta(r-)$. So the kernel of $\phi$ restricted to $R$ is precisely 
$$\{r \in R \, |\, \delta(rx) \in m, \, \forall x \in R, \, \forall \delta \in D_k(R)\}.$$
 For an element $x$ of $R$ and a $k$-linear differential operator $\delta$, the composition $\delta(x-)$ is also another differential operator. So the kernel of $\phi$ restricted to $R$ is  $\mathcal{M}_R$; see \Cref{l: larget D-ideal}. Since $\phi$, being $\widehat R$-linear, is obtained from $\phi_{|R}$ by extension of scalars to $\widehat R$, we conclude that the kernel of $\phi$ is $\mathcal{M}_R\hat R$ using  Lemma \ref{l: kernel of completion}.
\end{proof}

\begin{remark} \label{graded} 
To adapt the proof of Theorem \ref{c: main theorem} to the case of a  finitely generated graded $k$-algebra $(R, m)$ (again, with $k\cong R/m$), observe that in this case, the ring of differential operators $D_k(R)$ has a natural $\mathbb Z$-grading inherited from the natural grading on $\Hom_k^{cts}(R, R)$ 
\[
[\Hom_k^{cts}(R, R)]_{\ell} := \{R\xrightarrow{\phi} R \,\, | \,\, \phi(x)\in R_{\deg x + \ell} \,\,\, {\text{for all homogeneous }} x\in R\}.
\]
It follows that the ideal 
\[ \mathcal M_R :=\{x \in R \,\, | \,\, \phi(x)\in m\,\,\,  \forall \, \phi\in D_k(R)\}
\]
of  Proposition \ref{l: larget D-ideal}
is homogeneous, and   the maximal proper  $D_k(R)$-submodule of $R$. Likewise, 
$\Hom_k^{cts}(R, k)$ has a natural $\mathbb Z$-grading, and so because differential operators play nicely with localization (see \S \ref{localize}),
  the images of 
 $D_k(R)$ and  $D_k(R_m)$ in $D_k(R, k)=\Hom_k^{cts}(R, k)$ are the same.
Consequently,  the graded ring $(R, m)$ is a simple $D_k(R)$-module
 if and only if its localization $R_m$ is a simple $D_k(R_m)$-module, and  a graded version of Theorem  \ref{c: main theorem} holds.   In particular,  a graded version of Theorem \ref{MainThm} holds as well; see also \cite[Prop 3.2]{Jeff21}.

 \end{remark}
 
 \medskip
\begin{cor}\label{p: D-simplicity does not see completion}
Let  $(R,m)$ be a  Noetherian  local $k$-algebra such that $k\hookrightarrow R\to R/m$ is an isomorphism, and let $\widehat R$ denote its $m$-adic completion.
If $R$ is a simple $D_k(R)$-module, then  $\widehat R$ is a simple $D_k(\widehat R)$-module.
\end{cor}

\begin{proof}
Since $D_k(R) \subseteq D_k(\widehat R)$,  and $D_k(R, k) = E_R(k) = E_{\widehat R}(k) =  D_k(\widehat R,k)$, 
the surjectivity of $D_k(R) \rightarrow D_k(R, k) =E$
 implies surjectivity of $D_k(\widehat R) \rightarrow D_k(\widehat R,k)=E$.
 \end{proof}

The converse of Corollary \ref{p: D-simplicity does not see completion} 
does not  hold in general; We construct  examples in the next section. On the other hand, the converse {\it does hold} for sufficiently nice rings:

\begin{cor}\label{c: largest D ideal and completion}
Let  $(R,m)$ be a  Noetherian   local $k$-algebra such that $k\hookrightarrow R\to R/m$ is an isomorphism, and let $\widehat R$ denote its $m$-adic completion.
Assume, furthermore, that $R$ is 
essentially finite type over $k$ or that $R$ is $F$-finite.
 Then 
using  the notation from  Proposition \ref{l: larget D-ideal}, 
$$\mathcal M_{\widehat R}= \mathcal M_R \widehat R.$$ 
In particular, $R$ is simple as a $D_k(R)$ module if and only if $\widehat R$ is simple as a $D_k(\widehat R)$-module.
\end{cor}

\begin{proof}
Both $R$ and $\widehat R$ satisfy the hypothesis of Theorem \ref{MainThm}.  So the Matlis duals of the natural maps  
\begin{equation}\label{5}
D_k(R) \to D_k(R, k) \,\,\,\, \,\,\,\,\,{\text{and}}\,\,\,\,\,\,\,\, D_k(\widehat R) \to D_k(\widehat R, k)
\end{equation}
have kernels 
\[
 \mathcal M_R \widehat R\,\,\, \,\,\,\, \,\,\,\,\,{\text{and}}\,\,\,\,\,\,\,\, \,\,\,  \mathcal M_ {\widehat R},
\]
respectively. Note that the maps in (\ref{5}) have the same  target $$D_k(R, k) = \Hom_k^{cts}(R, k)=\Hom_k^{cts}(\widehat R, k) = D_k(\widehat R, k)  ,$$ which is the injective hull of the residue field $E$ for both $R$ and $\widehat R$. 

On the other hand, we have an isomorphism  $\widehat R\otimes_R D_k(R)\to D_k(\widehat R)$ by  \Cref{p: differential operator of a completion} in the finite type case or by Proposition \ref{p: for an F-finite ring differential operators of completion is the scalar extension of differential operators} in the $F$-finite case. 
So the second map in (\ref{5}) is canonically extended from the first:
\begin{equation}\label{7}
D_k(\widehat R) \cong \widehat R\otimes_R D_k(R) \to D_k(R, k)\cong D_k(\widehat R, k).
\end{equation}
This implies that the Matlis duals of the maps in (\ref{5})  can be viewed as   the {\it same } map
\[
\widehat R \to \Hom_R(D_k(R), E)=\Hom_{\widehat R}(\widehat R_R\otimes D_k(R), E) =  \Hom_{\widehat R}(D_k(\widehat R), E) 
\]
after making the natural identifications  of (\ref{9}) and  (\ref{7}). Hence, the  Matlis duals of the maps in (\ref{5})   have the same kernels--- that is, $\mathcal M_{\widehat R}= \mathcal M_R \widehat R.$

It follows,  in our setting,  that $R$ is $D$-simple if and only if $\widehat  R$ is $D$-simple, as $\mathcal M_R$ will be zero if and only if $\mathcal M_{\widehat R}$ is zero.
\end{proof}

\medskip


\section{Rings with few differential operators}

We now construct our cautionary examples of regular rings with few differential operators.
These examples are refinements of the discrete valuation rings constructed in \cite{DS18}. The construction generalizes to dimension two,  and possibly higher dimension;  see  Remark \ref{higherdim}.

\begin{construct}\label{example}
Fix any field $k$.  Choose  any power series $g \in k[[t]]$  such that $g(0)=0$ and $\frac{d}{dt}g\not\in k(t, g)$ (such $g$  always exist by  Proposition \ref{p: construction of g} below). Note that the field $L = k(t, g(t))$ is a subfield of fraction field $K = k[[t]][ \frac{1}{t}]$ of  $k[[t]]$. Restricting the $t$-adic valuation on $K$ to $L$, we get a   valuation on $L$ whose valuation ring\footnote{See \cite[Ch VI]{Bourbaki} or \cite[Ch  4]{Matsumura} for basics on valuation theory.} is 
 \begin{equation}\label{V}
 V_g= k(t,g) \cap k[[t]].
 \end{equation}
 We denote the maximal ideal of $V_g$ by  $\frm$.
  \end{construct}
 The  rings  $V_g$ constructed above have many nice properties:
 \begin{prop}\label{nice}
Let  $V_g$ be a ring as defined  in Construction \ref{example}. Then $(V_g, \frm)$ is a regular local  $k$-algebra of dimension one such that 
$k\hookrightarrow V_g \to V_g/\frm $ is an isomorphism. Furthermore
\begin{enumerate}
\item[(a)] The completion  $\widehat V_g$ of $V_g$ at its maximal ideal is isomorphic to $k[[t]]$.
\item[(b)] The  fraction field of $V_g$ is purely transcendental of transcendence degree two over $k$.
\end{enumerate}
 \end{prop}
 
 \begin{proof}  Since $V_g$ is the  valuation ring of a $\mathbb Z$-valued valuation on $L$, it is a discrete valuation ring in the field $L=k(t, g(t))$. In particular, $V_g$ is a one-dimensional regular local ring. Since the valuation is trivial on $k$, we see $k\subseteq V_g$, so that $V_g$ is  a $k$-algebra.  Note that the maximal ideal of  $V_g$ is 
$\frm = (t)\cap V_g$, and because the valuation takes the value $1$ on $t$, $\frm$ is generated by  $t$.
So
  the natural inclusions
 \[
 k\subseteq V_g \subseteq k[[t]]
 \]
  induce  isomorphisms
  \[
  k \hookrightarrow V_g/\frm \rightarrow k[[t]]/(t) \cong k.
  \]
Likewise, the inclusions 
 \[
 k[t]\subseteq V_g \subseteq k[[t]]
 \]
induce isomorphism of  $t$-adic completions
 \[
 k[[t]]\subseteq \widehat {V_g} \subseteq k[[t]],
 \]
 establishing (a). For  (b), it suffices to show that   the power series $g(t)$ is transcendental over $k(t)$. This  is a consequence of the 
 condition $\frac{dg}{dt}\not\in k(t,  g)$. Indeed, 
if there were an equation of integral dependence for $g$ over $k(t)$,  pick the minimal degree one
 \[
 g^n + a_{n-1}g^{n-1} + \cdots + a_1 g + a_0=0\,\,\,\,\,\,\,\,{\text{where  }} a_i\in k(t).
 \]
  Applying the $k$-linear derivation  $\frac{d}{dt}$ on  $K=k[[t]][ \frac{1}{t}]$,  we get an equation in $K$
 \[
 ng^{n-1}\frac{dg}{dt}  + (n-1)a_{n-1}g^{n-2} \frac{dg}{dt}  + \cdots + a_1 \frac{dg}{dt} + 
  \frac{da_{n-1}} {dt} 
  g^{n-1} + \cdots +  \frac{da_{1}}{dt} g + \frac{da_{0}}{dt} =0.
 \]
 Now note that $k(t)\hookrightarrow k((t)) $ is separable (because the ring $k[t]$ is excellent and hence its formal fibers are geometrically regular), so that $ng^{n-1} + (n-1)a_{n-1}g^{n-2}  + \cdots + a_1\neq 0$. 
 Factoring out $\frac{dg}{dt}$, we have 
 \begin{equation}\label{g}
 \frac{dg}{dt}=  
  \frac{ -(
a_{n-1}'
  g^{n-1} + \cdots + a_1'g + a_0')
   }
  {
  ng^{n-1} + (n-1)a_{n-1}g^{n-2}  + \cdots + a_1
  }  \in k(t, g),
 \end{equation}
 contrary to our assumption on  $g$. 
  So  the fraction field $L=k(t, g)$ of $V_g$ is a  rational function field of transcendence degree two over $k$.
 \end{proof}

The  $k$-algebras $V_g$ are pathological  from  the point of  view of differential operators:
    
\begin{theorem}\label{p: condition on derivative restricts differential operators}
 Fix any field  $k$. Let $g \in k[[t]]$ be a power series such that $g(0)=0$ and $\frac{d}{dt}g\not\in k(t, g)$. Then for  the discrete valuation ring 
  $V= k(t,g) \cap k[[t]]$,  the inclusion 
  \[V \rightarrow D_k(V)\]
   is an isomorphism.
\end{theorem}

The following corollary is immediate:
    
  \begin{cor} \label{pathological} Let $V_g$ be a regular local $k$-algebra as constructed in Construction \ref{example}. Then 
   $V_g$  is not simple  as module over  $D_k(V_g)$; indeed, it has infinite length over   $D_k(V_g)$.
  \end{cor}

  \begin{proof}[Proof of Corollary \ref{pathological}]
  Observe that any chain of ideals in $V_g$  (such as $\frm\supsetneq \frm^2 \supsetneq \frm^3\supsetneq \cdots)$  is a chain of $D_k(V_g)$-modules, since $D_k(V_g)=V_g$.     \end{proof}

Since nonzero derivations are differential operators of order one and not of order zero, the next corollary is immediate. 
\begin{cor}\label{c: DVRs with no derivations}
    Let $V_g$ be the discrete valuation $k$-algebra as constructed in Construction \ref{example}. Then $V$ does not admit any nontrivial $k$-linear derivation. 
\end{cor}  

\begin{remark}
Although the discrete valuation ring $V_g$ in the corollary above does not have any $k$-linear derivations, the module of K\"ahler differentials $\Omega_k(V)$ is nonzero. Indeed $\Omega_{V_g/k} \otimes_k L$- where $L= k(t, g(t))$ is the fraction field of $V_g$- is a two dimensional vector space over $L$. In fact $\Omega_{V_g/k}$ cannot be a finitely generated $V_g$-module, as $\textup{Hom}_{V_g}(\Omega_{V_g/k}, V_g) \otimes_{V_g} L \cong \textup{Der}_k(V_g) \otimes_{V_g} L $ is not isomorphic to  $\textup{Hom}_{L}(\Omega_{L/k}, L)$. Since $V_g$ linear maps from $\Omega_{V_g/k}$ to $V_g$ factors though $\Omega_{V_g/k}/(\Omega_{V_g/k})_{\text{tor}}$, the latter quotient module in fact is not finitely generated.
\end{remark}

\begin{proof}[Proof of Theorem \ref{p: condition on derivative restricts differential operators}] First observe that for every $n$, there is a canonical injection
 \begin{equation}\label{incl}
 D^n_k(V) \subseteq D^n_k(\widehat V)\subseteq D_k^n(K),
 \end{equation}
 where $K=k[[t]][\frac{1}{t}]$ is the fraction field of $\widehat V$.
  Indeed, the first  inclusion holds because each order $n$ differential  operator $V\to V$  extends to the operator   $V\to \widehat V$ in $D^n_k(V, \widehat V)$, and thus to a  unique operator in $D^n_k(\widehat V)$ by Proposition \ref{p: differential operator of completion without any noetherian assumption}; the second inclusion holds because  differential operators naturally extend to differential operators after localization (see \S \ref{localize}). 
  
We claim that
\begin{equation}\label{e: description of differential opeators of V}
    D^n_k(V)= \{\delta \in D^n_k(\widehat V) \, \, | \,\,  \delta(L) \subseteq L\},
\end{equation}
where $L=k(t, g)$ is the fraction  field of $V$, a subfield of $K$.
To see this, take  $\delta$ in the right hand side of (\ref{e: description of differential opeators of V}). Then  $\delta \in D^n_k(V)$ because $\delta(V)\subseteq  L\cap \widehat V = V$, and 
for $v_1, v_2, \ldots, v_{n+1} \in V$, $[\ldots[[\delta, v_1], v_2], \ldots, v_{n+1}]  =0$  as an operator on $\widehat  V$ and hence as an operator on $V$.  
The reverse inclusion is clear, as any differential operator on $V$ extends to a differential operator on its completion $\hat V$ by (\ref{incl}) and on its fraction field $L$ by (\ref{localization}). 
 This proves the claimed equality in (\ref{e: description of differential opeators of V}).\\

Next we claim that, if $V \subsetneq D_k(V)$, then $V \subsetneq D_k^1(V)$, so that $D_k(V)$  contains a differential operator of order precisely one. Indeed   supposing that $V \subsetneq D_k(V)$, let $n\geq  1$ be minimal such that $D^n_k(V) \setminus D^{n-1}_k(V)$ is non-empty.  
Choosing  $\delta \in D^n_k(V) \setminus D^{n-1}_k(V)$,  there must be  $f \in V$ such that $[\delta,f] \in D^{n-1}_k(V) \setminus D^{n-2}_k(V)$,  as otherwise   $\delta \in D^{n-1}_k(V)$. Thus  the minimal  $n=1$ and $V \subsetneq D_k^1(V)$.

We complete the proof of \Cref{p: condition on derivative restricts differential operators} by contradiction. Suppose $V \hookrightarrow D_k(V)$ is not an isomorphism. Then by the claim in the last paragraph,  we can  choose $\delta \in D^1_k(V) \setminus V$. Note that $D^1_k(V) \subseteq D^1_k(\widehat V) = D^1(k[[t]])$,  which is isomorphic to $k[[t]] \otimes_{k[t]}D^1(k[t]) $
by Corollary \ref{p: differential operator of a completion} applied to $S=k[t]$ and $J=(t)$. So since $D_k^1(k[t]) \cong k[t]\oplus k[t] \frac{d}{dt}$, every $k$-linear differential operator $\delta \in D_k(V)$  can be written 
$\delta= f_0 + f_1\frac{d}{dt}$ where $f_0, f_1 \in \widehat V. $  Since $\delta(1) \in V$, $f_0 \in V$. So also $\delta-f_0\in D^1_k(V)$, and 
thus $f_1= (\delta-f_0)(t) \in V$. Since $(\delta- f_0)(g)= f_1\frac{d}{dt}(g) \in V$, we conclude that $\frac{d}{dt}(g)$ is in the fraction field of $V$. This contradicts our assumption $\frac{dg}{dt} \notin k(t,g)$, and proves that $V\hookrightarrow D_k(V)$ is an isomorphism.\\
\end{proof}

Finally, we show that the $k$-algebras 
$V_g$  as constructed in Construction \ref{example} really exist:

\begin{prop}\label{p: construction of g}
 Let $k$ be a field, then there is a power series $g \in k[[t]]$ such that $g(0)=0$ and $\frac{dg}{dt}(g)\not\in k(t,g)$.    
\end{prop}

\begin{proof}
 First we treat  the case where $k$ has characteristic zero. Consider the power series
 $$g= \sum \limits_{j=0}^{\infty}(-1)^j\frac{t^{2j+1}}{(2j+1)!} \in \mathbb Q[[t]] ,$$
representing the entire function $\sin t$ on $\mathbb C$. We claim that $g$ is transcendental over $k(t)$.
 Indeed, $g$ is algebraic over $k(t)$ if and only if $k(t)[g]$ is a finite dimensional $k(t)$ vector space, and since $k(t)[g] \cong k \otimes_{\mathbb Q} \mathbb Q(t)[g]$, this is equivalent to saying that $ \mathbb Q(t)[g]$ has  finite dimension over $\mathbb Q(t)$. In turn, since 
  $\mathbb C\otimes_{\mathbb Q}  \mathbb Q(t)[g] \cong   \mathbb C(t)[g]$, 
this is equivalent to  $ \mathbb C(t)[g]  $ being finite dimensional  over  $ \mathbb C(t)$. Thus to check that $g$ is  transcendental over $k(t)$, it suffices to check $g$ is transcendental over $\mathbb C(t)$. This is clear: if $\sin t$ would be an algebraic function on $\mathbb C$, it could have only finitely many complex zeros, but we know that $\sin t$ has infinitely many.   Now observe that $\frac{dg}{dt} \notin k(t,g)$. Indeed, the formal derivative of the power series for $\sin t$ produces the power series for $\cos t$,  and in particular  these power series satisfy  $\sin^2 t+\cos^2 t =1$. Thus
$
\left(\frac{dg}{dt}\right)^2\in  k(t, g)
$
so that 
\begin{equation}\label{eqField}
k(t, g)\subseteq k(t,g,  \frac{dg}{dt} )
\end{equation}
 is a field extension of degree at  most  two. Note that $k(t, g)$ is the fraction field of a polynomial ring $A= k[t, g]$, since $g$ is transcendental over $k(t)$. So to prove that the extension (\ref{eqField}) is not of degree one, it suffices to show that the polynomial
 $X^2 + g^2-1$ is irreducible in the ring  $A[X]$. This is clear  by Eisenstein's criterion, using the prime ideal $(g-1)\subseteq  A=k[t, g]$ \cite[Exercise 17, \S 9.4]{DF}.

Now assume $k$ has characteristic $p>0$. 
We claim that $k[[t]]$ contains infinitely many elements algebraically independent over $k$. This   claim   reduces to the case where $k=\mathbb F_p$ as follows. Assuming that $\{x_i\}_{i\in I}$ is an infinite collection of algebraically independent power series in $\mathbb F_p[[t]]$, we have an embedding
\[
\mathbb F_p[\{x_i\}_{i\in I}] \hookrightarrow \mathbb F_p[[t]]
\]
of a 
polynomial ring in infinitely many variables into  $\mathbb F_p[[t]].$ Now applying the 
faithfully  flat base change $\mathbb F_p \hookrightarrow k$, we have
\[
k [\{x_i\}_{i\in I}] \hookrightarrow  k\otimes_{\mathbb F_p} \mathbb F_p[[t]] \subseteq k[[t]],
\]
which establishes that the $x_i$ are also algebraically independent when considered as power series over $k$.
Finally, to establish the case $k=\mathbb F_p$, observe that a field of finite transcendence degree over $\mathbb F_p$ is {\it countable}: 
indeed, the algebraic closure of a countable field is countable and any   {\it finitely generated}  transcendental extension of $\mathbb F_p$ is countable.  But $\mathbb F_p[[t]]$ is not countable, as its cardinality is the same as the set of infinite sequences of elements in $\mathbb F_p$.

Now, in light of the previous paragraph, let   $h, w \in k[[t]]$ be power series such that the set $\{t, w, h\}$  is algebraically independent over $k$, with $w(0)=0$.  Consider the power series  $g= th^p+ w^p$. Clearly $g(0)=0$. We claim  also that 
\[
\frac{d}{dt}(g) \notin k(t,g).\]
Indeed, if   $\frac{dg}{dt}= h^p \in k(t,g)$,  then also 
$w^p\in k(t,g)$. In particular, the algebraically independent set $\{t, h^p, w^p\} $ is contained in $k(t,g)$, a contradiction as the transcendence degree of $k(t,g)$ over $k$ is at most two.\footnote{In fact, the transcendence degree of $k(t, g)$ over $k$ is exactly two, by Proposition \ref{nice}.} This complete the proof of Proposition \ref{p: construction of g}.
\end{proof}

\begin{remark} \label{BadComplete} The $k$-algebra $V_g$ of Construction \ref{example} give examples of local  $k$-algebras  for which the ring of $k$-linear differential operators does not commute with completion. Indeed, since $D_k(V_g) = V_g$,  we see that 
\[
\widehat V_g \otimes_{V_g} D_k(V_g) \cong \widehat V_g,
\]
where $\widehat{\ }$  denotes the completion at the maximal ideal.
But 
\[
D_k(\widehat V_g) = D_k(k[[t]]),
\]
which contains  the derivation $\frac{d}{dt}$ (and many other differential operators of higher order on $\widehat V_g$).
\end{remark}

\begin{remark}\label{excZero}
The  $k$-algebras $V_g$ of Construction \ref{example} are excellent when $k$ has characteristic zero, they  are  Dedekind domains \cite[Cor 8.2.40]{Liu06}. However, if $k$ has prime characteristic, the ring $V_g$ is {\it not excellent}. Indeed, combining Remark \ref{BadComplete}   with
Proposition \ref{p: for an F-finite ring differential operators of completion is the scalar extension of differential operators}, we see that $V_g$ can not be $F$-finite, while for rings such as $V_g$, excellence is equivalent to $F$-finiteness \cite[Thm  2.4]{DS18}. 
\end{remark}

\begin{remark}\label{higherdim} Construction \ref{example} can be generalized to higher dimension as follows (with the caveat about Noetherianity   discussed below). Take a  power series $g\in k[[t]]]$ as in Proposition \ref{p: construction of g}. Consider the ring
\[
V_{d, g} =  k(t_1,  \dots, t_d, g(t_1), \dots, g(t_d)) \cap k[[t_1, \dots, t_d]]
\] 
obtained by intersecting two subrings inside the fraction field of $k[[t_1, \dots, t_d]].$  The arguments in Proposition  \ref{nice} show that $V_{d, g}$ is a  local $k$-algebra  with maximal ideal $\frm =(t_1, \dots, t_d)$  such  that $V_{d, g}/\frm \cong  k$, whose $\frm$-adic  completion is   $k[[t_1, \dots, t_d]]$ and whose fraction field  is  purely transcendental of  degree  $2d$    over $k$. 
However, $V_{d, g}$ is {\it not} always a valuation ring, so it is not obvious whether or not it is Noetherian. Our proof of Theorem \ref{p: condition on derivative restricts differential operators} does not require that $V$ be Noetherian, so it is also true that 
$D_k(V_{d,g})\cong  V_{d, g}$.  When $d=2$, the ring $V_{d, g}$ is Noetherian by \cite[Prop 3]{Valabrega}, so that $V_{2, g}$ is a regular local $k$ algebra of dimension two which is not simple as a $D_k(V_{d, g})$-module. More generally, the rings $V_{d,g}$ will be regular local $k$ algebras of dimension $d$ that are not $D$-simple if the following question has a positive answer:
\begin{question} \label{Q} Suppose that $g\in k[[t]]$ satisfies $\frac{dg}{dt}\notin k(t, g)$.  Is the ring
\[
V_{d, g} =  k(t_1,  \dots, t_d, g(t_1), \dots, g(t_d)) \cap k[[t_1, \dots, t_d]]
\] 
Noetherian for $d\geq 3$?
\end{question}
One might also experiment with taking {\it different} power series $g_i\in k[[t]]$ satisfying the conditions of Proposition \ref{p: construction of g}.
Related rings, including their  Noetherianity, have been studied by Heinzer, Rotthaus and Weigand,  but we believe that Question \ref{Q} is currently open. A positive answer to Question \ref{Q} for some $d$ will produce a regular local ring $V_{d, g}$ of dimension $d$ for which $D_k(V_{d, g})\cong  V_{d,g}$, so which is not $D$-simple.
\end{remark}

\section{Further Consequences in Prime Characteristic}

In prime characteristic, we get the  following related result:

\begin{theorem}\label{Indecomp}
   Fix a   perfect field $k$ of prime characteristic $p$. Take  any power series $g \in k[[t]]$  such that $g(0)=0$ and  $\frac{dg}{dt}\not\in k(t,g)$, and let   $V$ be the discrete valuation ring $k(t,g) \cap k[[t]]$. Then for each $n\in \mathbb N$, the inclusion $V \hookrightarrow {\End}_{V}(F^n_*V)$ is an isomorphism. \end{theorem}

\begin{remark}
 Theorem  \ref{Indecomp} does not follow from Theorem \ref{p: condition on derivative restricts differential operators} because    $D_k(R)$ need not be equal to $  \bigcup_{n\in\mathbb N} {\End}_{R}(F^n_*R)$  in general\footnote{ {\it A posteriori,}  however, for the non $F$-finite ring  $V$ it turns out that $D_k(V)\cong \bigcup_{n\in\mathbb N} {End}_{V}(F^n_*V)$, since both are just $V$.} for {\it non-$F$-finite} rings $R$. On the other hand, since the inclusion $D_k(R)\subseteq \bigcup_{n\in\mathbb N} {\End}_{R}(F^n_*R)$ does hold in general, Theorem \ref{Indecomp}  implies Theorem \ref{p: condition on derivative restricts differential operators}. 
\end{remark}

\begin{cor}
For the regular local ring $V$ as in Theorem \ref{Indecomp}, the modules $F_*^nV$ are indecomposable for all $n\in \mathbb N$.
\end{cor}

\begin{proof}[Proof of Corollary]
If $F_*^nV$ decomposed as $ M \oplus N$ for some nonzero submodules $M$ and $N$ of $F^n_*V$, then the projection $\pi$ onto, say $M$,  would be a homomorphism in ${\End}_{V}(F^n_*V)$ that is not simply multiplication by  some element $v\in V$. Indeed, $\pi$ kills the elements in $N$ but $F_*^nV$ is torsion free so $\pi$ can not be "multiplication by  $v\in V$" for any $v\in V$.
\end{proof}

\begin{proof}[Proof of Theorem \ref{Indecomp}]
Fix $n\in \mathbb N$. There is an natural inclusion 
\[
\End_V(F^n_*V) \hookrightarrow  \End_{\widehat V}(F_*^n\widehat V)
\]
obtained by applying the completion functor to each $\phi:F_*^nV\to F_*^nV$ (note that $\widehat {F_*^nV}$ is canonically identified with $F_*^n\widehat V$ by \cite[Lemma 1.14]{SchwedeSmith}).

Let $K$ be the fraction field of $\widehat V = k[[t]]$ and $L=k(t, g)$. Observe that 
\begin{equation}\label{e: description of Frobenius linear morphism of V}
  \text{End}_V(F^n_*V)= \{ \phi \in \text{End}_{K}(F^n_*K) \, \,| \, \, \phi (L)\subseteq L \, \, \text{and} \, \, \phi(\widehat V) \subseteq \widehat V \}.  
\end{equation}

 Since $k$ is perfect, the power series ring $\widehat   V = k[[t]]$ is $F$-finite, and hence so is  $K$. So $D_k(K)= \underset{n \in \mathbb{N}}{\cup}\text{End}_{K}(F^n_*K)$ by Proposition \ref{p: alternate description of the ring of differential operators of a prime characteristic ring}.
 So any $\phi \in \text{End}_{V}(F^n_*(V))$ is the restriction of a $k$-linear differential operator on $K$ sending $\widehat V$ to $\widehat V$ and $L$ to $L$. Hence any such $\phi$ is in $D_k(V)$, thanks to formula (\ref{e: description of differential opeators of V}). 
 Now by Theorem \ref{p: condition on derivative restricts differential operators}, we can conclude that  $\phi$ is multiplication by some element in $V$. So $V \hookrightarrow \text{End}_V(F^n_*V)$ is an isomorphism, completing the proof of Theorem \ref{Indecomp}.
\end{proof}

\bigskip
\begin{remark}
Theorem \ref{Indecomp} holds for the two-dimensional regular local rings $V_{2, g}$ as well, using the same proof; see Remark \ref{higherdim}. If Question \ref{Q} has a positive answer, then there are examples of regular local rings $R$  in every dimension such that $F_*^eR$ is indecomposable for every $e\in \mathbb N$.
\end{remark}

\begin{remark} If  $R$ is an {\it $F$-finite} regular local ring, such as a local ring of a point on a smooth variety of prime characteristic, 
the modules $F_*^nR$ are always {\it free} (of rank equal to $[K^{p^n}:K]$ where $K$ is the fraction field of $R$). So Theorem \ref{Indecomp} showcases the unexpected behavior of non-$F$-finite rings once again. 
For {\it non-regular} local rings, such as the local rings of {\it singular} points on algebraic varieties of prime characteristic, one might have  expected that "greater indecomposability" of  $F_*^nR$ can be viewed as a mark of "worse" singularities, and this appears to be the case for $F$-finite rings.

For example, 
consider a normal standard graded ring $R$ viewed as a the cone over a normal projective variety $X$. The $R$-modules $F_*^eR$ are $\frac{1}{p^e}$-graded and therefore always decompose\footnote{See \cite[\S 5.3]{SS}.} trivially as 
\begin{equation}\label{ST} 
F_*^eR = \bigoplus_{i=0}^{p^e-1} [F_*^eR]_{\frac{i}{p^e}\mod \mathbb Z},
\end{equation}
where $ [F_*^eR]_{\frac{i}{p^e}\mod \mathbb Z}$ denotes the $R$-submodule of  $F_*^eR$ consisting of elements whose degree is congruent to $\frac{i}{p^e}$ modulo $\mathbb Z$.
But 
for the cone $R$  over a smooth projective curve $X$ of 
genus at least two, each of the $p^e$ summands in (\ref{ST}) are indecomposable because the sheaves $F_*^n\mathcal L^i$ are stable  \cite[Prop 1.2]{LP} and hence indecomposable.  
 For a graded ring, this is  the "most indecomposable" we  can expect $F_*^nR$ to be.  Similarly, the singularity at the vertex of the cone gets more and more singular as the genus gets higher as well.
  For the same reason, the cone $R$ over an abelian variety with $p$-rank zero also has indecomposable $[F^n_*R]_{\mathbb Z} $  for all $n\in \mathbb N$ 
\cite[Thm 1.2]{ST}.  On the other hand, these examples still do not answer the following simple question:
\begin{question}
Does there exist a positive dimensional $F$-finite ring  $R$ such that $F^e_*R$ is indecomposable for all $e\in \mathbb N$?
\end{question}
Hochster and Yao ask further whether, for  any Noetherian module $M$ over a Noetherian $F$-finite ring,  the modules $F^e_*M$ are {\it decomposable}  for all $e\gg 0$. Hochster proved this in dimension one in \cite[Thm 5.16(2)]{Hochster}, and Yao (later with Hochster) showed the same in any dimension for rings $R$ that are (or have a finite extension that is) strongly $F$-regular \cite{Yao, HochYau}. See also \cite{SVdB}.
\end{remark}

\noindent{\bf Acknowledgements.} {We are grateful to Shunsuke Takagi  for directing our attention to the references \cite{LP} and \cite{ST},  to Sylvia Weigand and Christal Rotthaus for answering some questions about \cite{HRW}, and to Matthew Harrison-Trainer for clarifying some confusions about cardinalities. We thank Eamon Gallego and Linquan Ma for pointing towards useful references, Bhargav Bhatt and Jack Jeffries for their comments on an earlier draft. Jack Jeffries saved us from a misstatement in the proof of \Cref{c: main theorem}.}

\bibliographystyle{amsalpha}
\bibliography{Bibliography}
\end{document}